\documentclass[reqno]{amsart}
\usepackage{setspace}
\usepackage{graphicx}
\usepackage{amsmath}
\usepackage{amsfonts}
\usepackage{amssymb}

\begin{document}

{\large
\centerline{\Huge\bf Lower Bounds on }

\medskip
\centerline{\Huge\bf Lattice Covering Densities of Simplices}

\bigskip
\centerline{Miao Fu, Fei Xue and Chuanming Zong}

\vspace{0.6cm}
\centerline{\begin{minipage}{13cm}
{\bf Abstract.} This paper presents new lower bounds for the lattice covering densities of simplices by studying the Degree-Diameter Problem for abelian Cayley digraphs. In particular, it proves that the density of any lattice covering of a tetrahedron is at least $25/18$ and the density of any lattice covering of a four-dimensional simplex is at least $343/264$.
\end{minipage}}

\bigskip\noindent
\textbf{2020 Mathematics Subject Classification.} 52C17, 52B10, 52C07, 05C12.

\noindent
{\bf Keywords.} Lattice covering, simplex, the Degree-Diameter Problem.

\vspace{1.3cm}
\noindent
{\LARGE\bf 1. Introduction}

\medskip\noindent
More than 2,300 years ago, Aristotle claimed that {\it congruent regular tetrahedra can tile the whole space}. In other words, he believed that congruent regular tetrahedra can fill the whole space with neither gap nor overlap. Unfortunately, his claim was wrong: {\it Congruent regular tetrahedra can not tile the whole space}, which was discovered by Regiomontanus in the fifteenth century (see \cite{Laga12}). Then, two natural questions arose immediately: {\it What is the density of the densest tetrahedron packing and what is the density of the thinnest tetrahedron covering}? In fact, the packing case was emphasized by D. Hilbert \cite{Hilb00} as a part of his 18th problem.

Tetrahedron packings has been studied by many scholars, including mathematicians, physicists, and chemical engineers. Its history is dramatic and eventful. Let $S_3$ be a regular tetrahedron (a three-dimensional regular simplex). Let $\delta^c(S_3)$ and $\delta^l(S_3)$ denote the densities of the densest congruent packing and the densest lattice packing of $S_3$, respectively. Let's recall several events here, just as comparisons to the covering case. In 1904, H. Minkowski \cite{Mink04} claimed that $\delta^l(S_3)=\frac{9}{38}$. In fact, he also made a mistake. The correct answer is $\delta^l(S_3)=\frac{18}{49}$, which was proved by D.J. Hoylman \cite{Hoyl70} in 1970. Although Aristotle's claim was disproved by Regiomontanus in the fifteenth century which implied that $\delta^c(S_3)<1$, up to now we only know that $\delta^c(S_3)<1-2.6\times 10^{-25}$ which was discovered by S. Gravel, V. Elser and Y. Kallus \cite{Grav11} in 2011. For comprehensive surveys on packings and tetrahedron packings, we refer to \cite{Feje93} and \cite{Laga12}, respectively.

Covering is often regarded as a counterpart of packing. At least, their concepts, problems and results are often pairwise indeed. Nevertheless, covering came to mathematics much later than packing and our knowledge about it is still very limited. Let $K$ denote a convex body in $\mathbb{E}^n$ and let $C$ denote a centrally symmetric one. In particular, let $B_n$ denote the $n$-dimensional unit ball centered at the origin and let $S_n$ denote an $n$-dimensional regular simplex with unit edges and centered at the origin. Let $\theta^c(K)$, $\theta^t(K)$ and $\theta^l(K)$ denote the densities of the thinnest congruent covering, the thinnest translative covering, and the thinnest lattice covering of $\mathbb{E}^n$ with $K$, respectively. Clearly, for every convex body $K$ we have
$$1\leq \theta^c(K)\leq \theta^t(K)\leq \theta^l(K). \nonumber $$
Moreover, both $\theta^t(K)$ and $\theta^l(K)$ are invariant under non-singular linear transformations on $K$. Therefore, in this paper we will work on the simplex with vertices $(0,0, \ldots, 0)$, $(1, 0, \ldots, 0)$, $(0,1, \ldots, 0)$, $\ldots$, $(0, 0, \ldots, 1)$, instead of the regular one.

In the plane, our covering knowledge is comparatively complete. In 1939, R. Kershner \cite{Kers39} proved that
$$\theta^t(B_2)=\theta^l(B_2)= \frac{2\pi}{\sqrt{27}}.$$
In 1946 and 1950, L. Fejes T\'oth \cite{Feje46,Feje50} showed that
$$\theta^t(C)=\theta^l(C)\le \frac{2\pi}{\sqrt{27}}$$
holds for all centrally symmetric convex domains $C$ and the second equality holds if and only if $C$ is an ellipse.
In 1950, I. F\'{a}ry \cite{Fary50} proved that
$$\theta^l(K)\le \frac{3}{2}$$
holds for all convex domains $K$, where the equality holds if and only if $K$ is a triangle. It is rather surprising that
$$\theta^t(S_2)=\frac{3}{2}$$
was proved only in 2010 by J. Januszewski \cite{Janu10}. It is even more surprising that, up to now some basic covering problems in the plane are still open (see \cite{Zong14}). For example, we do not know yet if $\theta^t(K)=\theta^l(K)$ holds for all convex domains.

In $\mathbb{E}^3$, we only know one exact covering result
$$\theta^l(B_3)= \frac{5\sqrt5\pi}{24},$$
which was discovered in 1954 by R. P. Bambah \cite{Bamb54}. For tetrahedron coverings, several bounds have been achieved.
In the lattice case, we have
$$\frac{2^{16}+1}{2^{16}}\le \theta^l(S_3)\le \frac{125}{63},$$
where the upper bound was discovered by C. M. Fiduccia, R. W. Forcade and J. S. Zito \cite{Fidu98} and R. Dougherty and V. Faber \cite{Doug04} in the 1990s by constructing a particular lattice covering and the lower bound was achieved by F. Xue and C. Zong \cite{Xue18} in 2018. In fact, it was proved by R. Forcade and J. Lamoreaux \cite{Forc00} and R. Dougherty and V. Faber \cite{Doug04} that the upper bound is a local minimum. They even conjectured it to be the exact value of $\theta^l(S_3)$. In the congruent case, in 2006 J. H. Conway and S. Torquato \cite{Conw06} obtained
$$\theta^c(S_3)\le \frac{9}{8}$$
by constructing a particular tetrahedron covering. Nothing nontrivial is known about $\theta^t(S_3)$.

In $\mathbb{E}^n$, covering has been studied by R. P. Bambah, H. S. M. Coxeter, H. Davenport, P. Erd\H{o}s, L. Few, G. L. Watson, and in particular by C. A. Rogers (see \cite{Roge64}). They proved that
$$\theta^l(K)\leq n^{\log_2\log_e n+c},$$
$$\theta^t(K)\leq n\log n+n\log \log n+5n,$$
and
$$\frac{n}{e\sqrt{e}}\ll \theta^t(B_n)\leq \theta^l(B_n)\leq c\cdot n(\log_e n)^{\frac{1}{2}\log_2 2\pi e}.$$
Since the 1960s, progress in covering is very limited (see \cite{Bras05}). In 2018, F. Xue and C. Zong \cite{Xue18} discovered that
$$\theta^l(S_n)\geq 1+\frac{1}{2^{3n+7}}.\eqno(1.1)$$
In 2021, O. Ordentlich, O. Regev and B. Weiss \cite{Orde21} improved Rogers' upper bound to
$$\theta^l(K)\leq c n^2,$$
where $c$ is a suitable positive constant.

\medskip
In this paper, we prove the following results:

\medskip\noindent
{\bf Theorem 1.1.} {\it If $S_3+\Lambda$ is a lattice covering of $\mathbb{E}^3$, then its density is at least $25/18$. In other words, we have}
$$\theta^l(S_3)\ge \frac{25}{18}.$$

\smallskip\noindent
{\bf Theorem 1.2.} {\it If $S_4+\Lambda$ is a lattice covering of $\mathbb{E}^4$, then its density is at least $343/264$. In other words, we have}
$$\theta^l(S_4)\ge \frac{343}{264}.$$

\medskip
Our method is based on the {\it Degree-Diameter Problem} for {\it abelian Cayley digraphs}. Let $G$ be a finite {\it abelian group} and let $E=\{{\bf g}_1, {\bf g}_2, \ldots , {\bf g}_n\}$ be a set of generators for $G$. Then the abelian Cayley digraph of $G$ and $E$ has the elements of $G$ as its vertices and directed edges from each vertex ${\bf u}$ to all vertices ${\bf v}={\bf u}+{\bf g}_i$, $1\leq i\leq n$. Note that every vertex in the digraph has out-degree $n$, and the order of the digraph is $|G|$. Then, the Degree-Diameter Problem for abelian Cayley digraphs can be stated as:
{\it Given positive integers $n$ and $d$, find the largest order $f(n,d)$ among all abelian Cayley digraphs of $G$ and $E$, where $|E|=n$ and the diameters of these digraphs are at most $d$}.

Both the Degree-Diameter Problem and Cayley graphs are well-known in the graph theory community. For a survey and recent progress, we refer to M. Miller and J. \v{S}ir\'{a}\v{n} \cite{Mill05} and T. Zhang and G. Ge \cite{Zhang19}, respectively. In 1974, C. K. Wong and D. Coppersmith \cite{Wong74} showed that
$$f(2,d)=\left\lfloor\frac{(d+2)^2}{3}\right\rfloor. $$
In 1998, C. M. Fiduccia, R. W. Forcade and J. S. Zito \cite{Fidu98} proved that
$$f(3,d)\leq \frac{3(d+3)^3}{25}\eqno(1.2)$$
and asked for a similar bound for $f(4,d)$. In 2004, R. Dougherty and V. Faber \cite{Doug04} obtained that
$$\frac{cd^n}{n!n(\ln{n})^{1+\log_2e}}+O(d^{n-1})\leq f(n,d)\leq \dbinom{d+n}{n}.\eqno(1.3)$$

In this paper, we also prove the following result.

\medskip\noindent
{\bf Theorem 1.3.} {\it If an abelian Cayley digraph has degree at most $4$ and diameter at most $d$, then it has at most $\frac{11(d+4)^4}{343}$ elements. In other words,}
$$f(4,d)\leq \frac{11(d+4)^4}{343}.$$

\medskip\noindent
{\bf Remark 1.1.} Theorem 1.3 will be useful in the proof of Theorem 1.2. As one can check, Theorem 1.3 is better that (1.3) when $d$ is large. In fact, by similar method the upper bound in (1.3) can be improved to
$$f(n, d)\leq\frac{(d+n)^n}{n\cdot n!}\left(n-1+\left(\frac{n-1}{2n-1}\right)^{n-1}\right).\nonumber$$
Consequently, the lower bound in (1.1) can be improved to
$$\theta^l(S_n)\geq \frac{n}{n-1+\left(\frac{n-1}{2n-1}\right)^{n-1}}. \nonumber$$
However, since the improvements are not essential, we will not include their proofs here.

\vspace{1cm}
\noindent
{\LARGE\bf 2. The Degree-Diameter Problem}

\medskip\noindent
The Degree-Diameter Problem for abelian Cayley digraphs is closely related to lattice coverings (see \cite{Doug04} and \cite{Fidu98}). Let $\mathbb{Z}^n$ be the {\it integer lattice} generated by the standard basis $\{{\bf e}_1, {\bf e}_2, \ldots, {\bf e}_n\}$. From the algebraic point of view, $\mathbb{Z}^n$ is a free abelian group with $n$ generators. Thus for any finite abelian group $G$ generated by $\{{\bf g}_1, {\bf g}_2, \ldots, {\bf g}_n\}$ (in additive), there is a homomorphism $\phi: \mathbb{Z}^n\to G$ defined by
$$\phi\left( \sum_{i=1}^nz_i{\bf e}_i\right)= \sum_{i=1}^nz_i{\bf g}_i,\qquad z_i\in \mathbb{Z}.$$
Clearly, $\phi$ is surjective. Let $L_n$ be the {\it kernel} of $\phi$. It is well-known and easy to check that $L_n$ is an $n$-dimensional sublattice of $\mathbb{Z}^n$. Then there is an isomorphism
$$\overline{\phi}: \mathbb{Z}^n/L_n\to G $$
and
$$|G|=\left[\mathbb{Z}^n\! :L_n\right]= d(L_n),$$
where $d(L_n)$ denote the determinant of $L_n$.

Let ${\bf o}$ be the origin of $\mathbb{E}^n$ and let $\mathcal{O}_n$ denote the positive orthant in $\mathbb{E}^n$. We start from ${\bf o}$ and perform a breadth-first search in $\mathbb{Z}^n\cap\mathcal{O}_n$ until we find all representative elements of the group $G$. Let ${\bf x}=(x_1, x_2,\dots, x_n)$ and ${\bf y}=(y_1, y_2,\dots, y_n)$ be two distinct points in $\mathbb{Z}^n\cap\mathcal{O}_n$ and let $|{\bf o}, {\bf x}|_M$ denote the {\it Manhattan distance} between ${\bf o}$ and ${\bf x}$. In other words,
$$|{\bf o}, {\bf x}|_M=\sum_{i=1}^nx_i.$$
Then, we write ${\bf x}\prec {\bf y}$ if either $|{\bf o}, {\bf x}|_M<|{\bf o}, {\bf y}|_M$ or $|{\bf o}, {\bf x}|_M=|{\bf o}, {\bf y}|_M$ but ${\bf x}$ precedes ${\bf y}$ lexicographically. For example, in $\mathbb{Z}^2\cap\mathcal{O}_2$ we have
$$(0,0)\prec (0,1)\prec (1,0)\prec (0,2)\prec (1,1)\prec (2,0)\prec (0,3)\prec (1,2)\prec (2,1)\prec (3,0)\prec\ldots $$
and in $\mathbb{Z}^3\cap\mathcal{O}_3$ we have
$$(0,0,0)\prec (0,0,1)\prec (0,1,0)\prec (1,0,0)\prec (0,0,2)\prec (0,1,1)\prec (0,2,0)\prec (1,0,1)\prec \ldots $$
Clearly, the linear order $\prec $ compatible with addition. In other words, if ${\bf x}\prec {\bf y}$, then
$${\bf x}+{\bf z}\prec {\bf y}+{\bf z}$$
holds for all ${\bf z}\in \mathbb{Z}^n\cap\mathcal{O}_n$.

For ${\bf g}\in G$, let $\phi^*({\bf g})$ be the first lattice point ${\bf x}$ in the ordered sequence satisfying $\phi({\bf x})={\bf g}$. Then
$$T_n=\{\phi^*({\bf g}): {\bf g}\in G\}\eqno(2.1)$$
is a complete set of coset representatives for the kernel $L_n$. From the geometric point of view, $T_n+L_n$ is a {\it tiling} of $\mathbb{Z}^n$, i.e., $T_n+L_n=\mathbb{Z}^n$ and the translates are pairwise disjoint. Of course, we have
$$|G|=\left[\mathbb{Z}^n\! :L_n\right]=|T_n|.\eqno(2.2)$$
Note that the diameter of the abelian Cayley digraph of $G$ and $E=\{{\bf g}_1, {\bf g}_2, \ldots, {\bf g}_n\}$ is equal to the greatest Manhattan distance of any point in $T_n$ from ${\bf o}$. We will accordingly call it the {\it M-diameter} of $T_n$. Thus we have the following result:

\medskip\noindent
{\bf Proposition 2.1.} {\it The abelian Cayley digraph of $G$ and $E$ has a diameter of at most $d$ if and only if $T_n+L_n$ is a tiling of $\mathbb{Z}^n$, where $T_n$ has a M-diameter of at most $d$. Moreover, the order of the graph is equal to $|T_n|$, i.e., $\left[\mathbb{Z}^n\! :L_n\right]$. }

\medskip
If $T_n+L_n$ is a tiling of $\mathbb{Z}^n$, then $\overline{T_n}+L_n$ is a tiling of $\mathbb{E}^n$, where
$$\overline{T_n}=T_n+[0, 1)^n \nonumber $$
is called a {\it Cayley tile} of $\mathbb{E}^n$. Clearly, we have
$${\rm vol}(\overline{T_n})=|T_n|=d(L_n),\eqno(2.3)$$
where ${\rm vol}(\overline{T_n})$ denotes the volume of $\overline{T_n}$. If $d$ is the M-diameter of $T_n$, then the greatest Manhattan distance between a point of $\overline{T_n}$ and ${\bf o}$ is $d+n$. We say that the {\it M-diameter} of $\overline{T_n}$ is $d+n$. Based on Proposition 2.1 and (2.3), the order of the graph is equal to ${\rm vol}(\overline{T_n})$, and the Degree-Diameter Problem for abelian Cayley digraphs is equivalent to the following problem:
\begin{quote}
{\it Given positive integers $n$ and $d$, find the largest ${\rm vol}(\overline{T_n})$ among all Cayley tiles $\overline{T_n}$ of $\mathbb{E}^n$ with $M$-diameter at most $d+n$.}
\end{quote}

Let's recall some useful definitions and notions introduced by Fiduccia, Forcade and Zito \cite{Fidu98}. Assume that both ${\bf x}=(x_{1}, x_{2},\dots, x_{n})$ and ${\bf y}=(y_{1}, y_{2},\dots, y_{n})$ are points in $\mathcal{O}_n$. We say that ${\bf x}\preceq {\bf y}$ if $x_{i}\leq y_{i}$ holds for all $i$. A {\it notch} in the boundary of $\overline{T_n}$ is a place where it looks like a translation of $\mathcal{O}_n$ has been cut out of the tile (for example, the point ${\bf v}$ in Fig. 1). The {\it silhouette} of $\overline{T_n}$ is the set of points ${\bf p}$ with at least one zero coordinate such that ${\bf p}\preceq {\bf q}$ holds for some point ${\bf q}\in \overline{T_n}$. In fact, as one can see from the next lemma, it is the union of the projections of $\overline{T_n}$ onto the coordinate hyperplanes. Let $\pi_i$ denote the projection from $\mathbb{E}^n$ to the hyperplane $H_i=\{{\bf x}\! :\ x_i=0\}$. Then, as shown in Fig. 1, the silhouette of $\overline{T_n}$ is $\bigcup_{i=1}^{n} \pi_i(\overline{T_n})$.

\begin{figure}
\centering
\includegraphics[height=4.5cm]{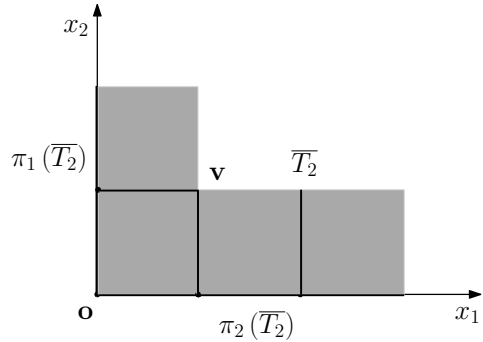}
\caption{A Cayley tile in $\mathbb{E}^2$ and its silhouette}
\end{figure}

\medskip\noindent
{\bf Lemma 2.1 (Fiduccia, Forcade and Zito \cite{Fidu98}).} {\it The Cayley tile $\overline{T_n}$ has the following properties:
\begin{enumerate}
\item If ${\bf x}\in \overline{T_n}$ and ${\bf o}\preceq {\bf y}\preceq {\bf x}$, then ${\bf y}\in \overline{T_n}$.
\item It has at most one notch.
\item It is uniquely determined by its silhouette and, if it has a notch, the coordinates of the notch.
\end{enumerate}}

\medskip\noindent
{\bf Remark 2.1.} Let ${\bf v}$ be a notch of $\overline{T_n}$. The last properties of Lemma 2.1 implies that, if all projections $\pi_i({\bf p})$ of a point ${\bf p}$ are in the silhouette $\overline{T_n}$ and ${\bf v}\not \preceq {\bf p}$, then ${\bf p}\in \overline{T_n}$. Let $[{\bf p},\pi_i({\bf p})]$ denote the segment connecting ${\bf p}$ and $\pi_i({\bf p})$. It follows from the first property of Lemma 2.1 that $\pi_i({\bf p})$ belongs to the silhouette of $\overline{T_n}$ if and only if
$$[{\bf p},\pi_i({\bf p})]\cap \overline{T_n}\not=\varnothing .$$

\medskip
Recall that $\overline{T_4}$ is a Cayley tile of $\mathbb{E}^4$ with $M$-diameter of $d^*$, where $d^*\le d+4$. We proceed to prove Theorem 1.3 by dealing with two cases:

\medskip\noindent
{\bf Case 1. The Cayley tile $\overline{T_4}$ has no notch}

\medskip
By Lemma 2.1, $\overline{T_4}$ is completely determined by its silhouette. For each point ${\bf p}\notin \overline{T_4}$ in $\mathcal{O}_{4}$, according to Remark 2.1, there must be a $\pi_i$ satisfying
$$[{\bf p},\pi_i({\bf p})]\cap \overline{T_4}=\varnothing.$$
For convenience, then we say ${\bf p}$ is of $\pi_i$-{\it type}. Note that there are at most four types and a point can be of more than one type simultaneously.

\medskip\noindent
{\bf Lemma 2.2.} {\it Assume that ${\bf p}$ and ${\bf q}$ are two points in $\mathcal{O}_{4}$. If ${\bf p}\preceq {\bf q}$ and ${\bf p}$ is of $\pi_i$-type, then ${\bf q}$ is also of $\pi_i$-type.}

\medskip\noindent
{\bf Proof.} Suppose that ${\bf p}$ is $\pi_i$-type. Then, we have
$$[{\bf p},\pi_i({\bf p})]\cap \overline{T_4}=\varnothing.$$
Therefore, we have
$$\pi_i({\bf p})\notin \overline{T_4}.$$
Since ${\bf p}\preceq {\bf q}$, for every point ${\bf w}\in [{\bf q},\pi_i({\bf q})]$ one can deduce that
$$\pi_i({\bf p})\preceq {\bf w}.$$
Thus ${\bf w}$ can not be a point of $\overline{T_4}$, which implies that
$$[{\bf q},\pi_i({\bf q})]\cap \overline{T_4}=\varnothing.$$
Therefore, ${\bf q}$ is also $\pi_i$-type. The lemma is proved. \hfill{$\Box$}

\bigskip
For convenience, we define
$$S_{4,d^*}=\{(x_{1}, x_{2}, x_{3}, x_{4})\in \mathbb{E}^4\! :\ x_{i}\geq 0,\ x_{1}+x_{2}+x_{3}+x_{4}\leq d^*\}. $$
Note that the greatest Manhattan distance between any point of $S_{4,d^*}$ and ${\bf o}$ is $d^*$, which implies that
$$\overline{T_4}\subseteq S_{4,d^*}.\eqno(2.4)$$
Clearly,
$$F=\{(x_{1}, x_{2}, x_{3}, x_{4})\in \mathbb{E}^4\! :\ x_{i}\geq 0,\ x_{1}+x_{2}+x_{3}+x_{4}= d^*\}$$ is a facet of $S_{4,d^*}$, which in fact is a regular tetrahedron.

Let $F_1$ denote the subset of $F$ including all points of $\pi_1$-type and successively define $F_i$ to be the subset of $F$ including all points of $\pi_i$-type which are not in $F_1\cup \cdots \cup F_{i-1}$. Clearly, $F_1$, $F_2$, $F_3$ and $F_4$ are pairwise disjoint and, since $\overline{T_4}$ has no notch,
$$F=F_1\cup F_2\cup F_3\cup F_4.$$
For $1\le i\le 4$, we define
$$D_i=\bigcup_{{\bf p}\in F_i}[{\bf p},\pi_i({\bf p})].$$
Since $[{\bf p}, \pi_i({\bf p})]\cap \overline{T_4}=\varnothing$ holds for all ${\bf p}\in F_i$, we have
$$D_i\cap \overline{T_4}=\varnothing,$$
which implies that $D_i$ is in the complement of $\overline{T_4}$.

\medskip\noindent
{\bf Lemma 2.3.} {\it The sets $D_{i}$, where $1\leq i\leq 4$, are pairwise disjoint.}

\medskip\noindent
{\bf Proof.} If, on the contrary, $D_i$ and $D_j$ have a common point ${\bf p}$, where $i<j$. Then there is a point ${\bf p}_i\in F_i$ and a point ${\bf p}_j\in F_j$ such that
$${\bf p}\in [{\bf p}_i, \pi_i({\bf p}_i)]\cap [{\bf p}_j, \pi_j({\bf p}_j)],$$
where $[{\bf p}_i, \pi_i({\bf p}_i)]\cap \overline{T_4}=\varnothing$ and $[{\bf p}_j, \pi_j({\bf p}_j)]\cap \overline{T_4}=\varnothing$. Consequently, ${\bf p}\preceq {\bf p}_i$, ${\bf p}\preceq {\bf p}_j$, and ${\bf p}$ is both $\pi_i$-type and $\pi_j$-type.  According to Lemma 2.2, ${\bf p}_j$ is also of $\pi_i$-type, i.e., ${\bf p}_j\in F_i$, which contradicts the fact
$$F_i\cap F_j=\varnothing.$$
The lemma is proved. \hfill{$\Box$}

\medskip
By (2.4) and Lemma 2.3 we have
$${\rm vol}(\overline{T_4})\leq {\rm vol}(S_{4,d^*})-\sum_{i=1}^{4}{\rm vol}(D_i)=\frac{d^{*4}}{4!}-\sum_{i=1}^{4}{\rm vol}(D_i).\eqno(2.5)$$
Notice that
$$\sum\limits_{i=1}^{4}{\rm vol}(D_{i})=\iiint\limits_F \lambda f({\bf x})d{\bf x}=\sum\limits_{i=1}^{4}\iiint\limits_{\pi_i(F_i)}\left( d^*\!-|{\bf o}, {\bf x}|_M\right) d{\bf x}, \eqno(2.6)$$
where $\lambda$ is a constant because we are integrating over $F$ instead of over the coordinate hyperplanes, $f({\bf x})$ is the distance from ${\bf x}$ to an appropriate coordinate hyperplane. In fact, if ${\bf x}=(x_1, x_2, x_3, x_4)\in F_i$, then
$$f({\bf x})=x_i.$$

To maximize the right-hand side of (2.5), we need to minimize the value of (2.6). For any ${\bf x}=(x_1, x_2, x_3, x_4)$ in $F$, we define
$$g({\bf x})=\min\!\ \{x_1, x_2, x_3, x_4\}.$$
Furthermore, we define
$$F'_i=\{{\bf x}\in F\! :\ g({\bf x})=x_i\}.$$
Of course, $\pi_{i}(F'_i)$ is the projection of $F'_{i}$ onto the nearest coordinate hyperplane. Thus, we have
\begin{align}
{\rm vol}(\overline{T_4})&\leq\frac{d^{*4}}{4!}-\iiint\limits_F \lambda g({\bf x})d {\bf x} \nonumber\\
&=\frac{d^{*4}}{24}-\sum_{i=1}^{4}\iiint\limits_{\pi_{i}(F'_i)} \left(d^*\!-|{\bf o}, {\bf x}|_M\right)d{\bf x} \nonumber\\
&=\frac{d^{*4}}{24}-4\iiint\limits_{\pi_{4}(F'_4)} \left(d^*\!-\!x_{1}\!-\!x_{2}\!-\!x_{3}\right) d{\bf x}. \nonumber
\end{align}
The integral region $\pi_{4}(F'_4)$ and its projection on the $x_{1}x_{2}$-coordinate plane are shown in Fig. \ref{4-1}. It follows that
%\begin{small}
\begin{align}
\iiint\limits_{\pi_{4}(F'_4)} \left(d^*\!-\!x_{1}\!-\!x_{2}\!-\!x_{3}\right) d{\bf x}
&=\int_{0}^{\frac{d^*}{4}}dx_{1}\int_{d^*\!-\!3x_{1}}^{d^*\!-\!x_{1}}dx_{2}\int_{\frac{1}{2}(d^*\!-\!x_{1}\!-\!x_{2})}^{d^*\!-\!x_{1}\!-\!x_{2}} \left(d^*\!-\!x_{1}\!-\!x_{2}\!-\!x_{3}\right) dx_{3} \nonumber\\
&\quad +\int_{\frac{d^*}{4}}^{d^*}dx_{1}\int_{\frac{1}{3}(d^*\!-\!x_{1})}^{d^*\!-\!x_{1}}dx_{2}\int_{\frac{1}{2}(d^*\!-\!x_{1}\!-\!x_{2})}^{d^*\!-\!x_{1}\!-\!x_{2}} \left(d^*\!-\!x_{1}\!-\!x_{2}\!-\!x_{3}\right) dx_{3} \nonumber\\
&\quad +\int_{0}^{\frac{d^*}{4}}dx_{1}\int_{x_{1}}^{d^*\!-\!3x_{1}}dx_{2}\int_{d^*\!-\!2x_{1}\!-\!x_{2}}^{d^*\!-\!x_{1}\!-\!x_{2}} \left(d^*\!-\!x_{1}\!-\!x_{2}\!-\!x_{3}\right) dx_{3} \nonumber\\
&\quad +\int_{0}^{\frac{d^*}{4}}dx_{2}\int_{x_{2}}^{d^*\!-\!3x_{2}}dx_{1}\int_{d^*\!-\!x_{1}\!-\!2x_{2}}^{d^*\!-\!x_{1}\!-\!x_{2}} \left(d^*\!-\!x_{1}\!-\!x_{2}\!-\!x_{3}\right) dx_{3} \nonumber\\
&=\frac{d^{*4}}{384}\nonumber.
\end{align}
%\end{small}

\begin{figure}
\centering
\includegraphics[height=5cm]{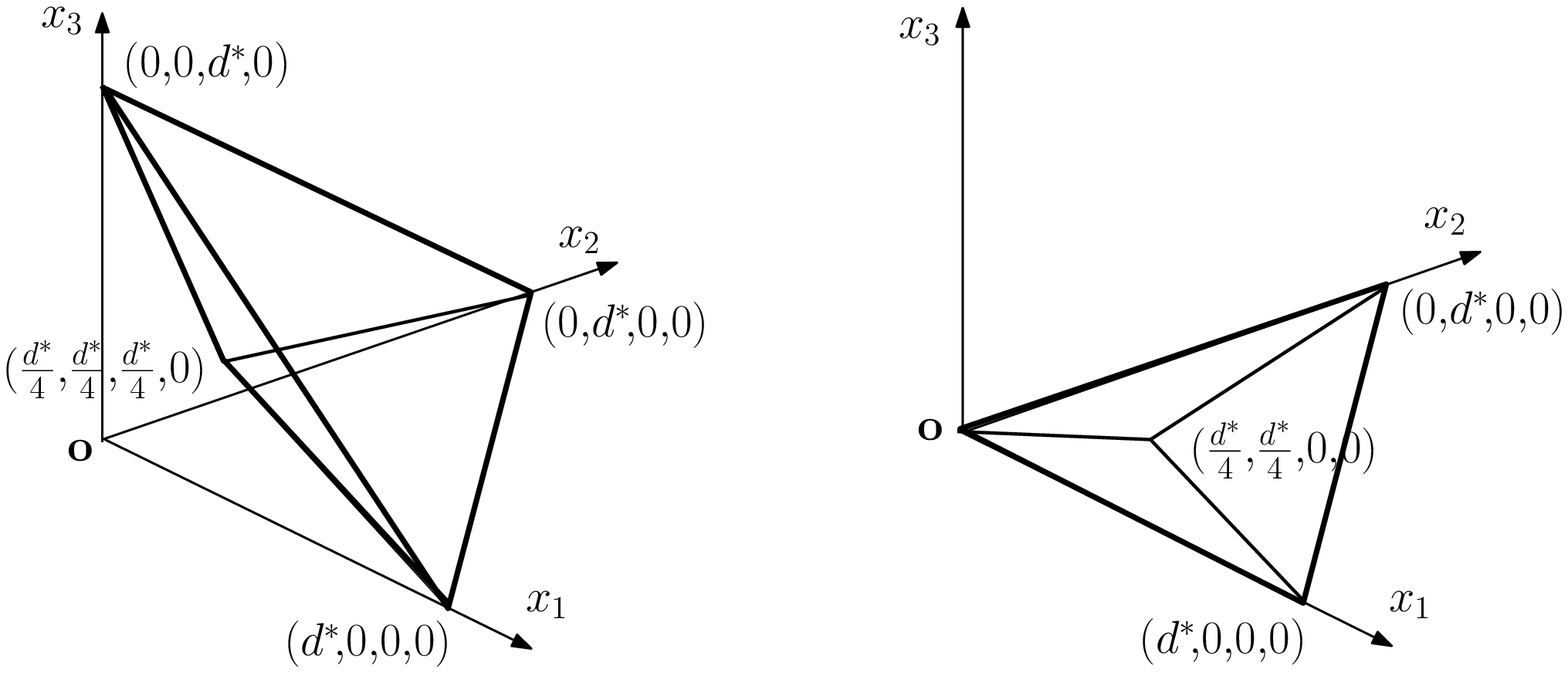}
\caption{The integral region $\pi_{4}(F'_4)$ and its projection on the $x_{1}x_{2}$-plane}
\label{4-1}
\end{figure}

Thus, with no notch in $\overline{T_4}$, we have
$${\rm vol}(\overline{T_4}) \leq\frac{d^{*4}}{32}. $$

\medskip\noindent
{\bf Case 2. The Cayley tile $\overline{T_4}$ has a notch}

\medskip
From Remark 2.1, if $\overline{T_4}$ has a notch ${\bf v}$, then not all points ${\bf x}\in F$ are of $\pi_1$-, $\pi_2$-, $\pi_3$-, and $\pi_4$-type, and those that are not must satisfy ${\bf v}\preceq {\bf x}$. Then, we define
$$Q=\{ {\bf x}\in F\! :\ {\bf v}\preceq {\bf x}\}.$$
Following the previous argument, let $F_1$ denote the subset of $F\setminus Q$ including all points of $\pi_1$-type and successively define $F_i$ to be the subset of $F\setminus Q$ including all points of $\pi_i$-type which are not in $F_1\cup \cdots \cup F_{i-1}$. Then, we write
$$D_i=\bigcup_{{\bf x}\in F_i}[{\bf x}, \pi_i({\bf x})].$$
It is easy to see that
$$F=F_1\cup F_2\cup F_3\cup F_4\cup Q$$
is a disjoint union, and the sets $D_i$ are in the complement of $\overline{T_4}$. Furthermore, similar to Lemma 2.3, one can prove that $D_i$ are pairwise disjoint. Then, we define
$$P=\left\{{\bf p}\in S_{4,d^*}\! :\ {\bf v}\preceq {\bf p}\right\}.$$
In other words, $P$ is the union of all segments from $Q$ to the notch point ${\bf v}$. Thus, we have
$${\rm vol}(\overline{T_4})\leq\frac{d^{*4}}{4!}-\sum\limits_{i=1}^{4}{\rm vol}(D_i)-{\rm vol}(P).\eqno(2.7)$$
Similar to (2.6), we have
$$\sum\limits_{i=1}^{4}{\rm vol}(D_{i})=\iiint\limits_{F\backslash Q} \lambda f({\bf x})d{\bf x}=\sum\limits_{i=1}^{4}\iiint\limits_{\pi_{i}(F_i)} \left(d^*\!-|{\bf o}, {\bf x}|_M\right) d{\bf x},\eqno(2.8)$$
where $\lambda$ is a constant because we are integrating over $F$ instead of over the coordinate hyperplanes, $f({\bf x})$ is the distance from ${\bf x}$ to an appropriate coordinate hyperplane. In fact, if ${\bf x}=(x_1, x_2, x_3, x_4)\in F_i$, then
$$f({\bf x})=x_i.$$

First, let's fix the volumes of $Q$ and $P$ and to minimize the value of (2.8). For a point ${\bf x}=(x_1, x_2, x_3, x_4)$ in $F\backslash Q$, we define
$$g({\bf x})=\min\!\ \{x_1, x_2, x_3, x_4\},$$
which is the distance from ${\bf x}$ to the nearest coordinate hyperplane, and
$$F'_i=\{{\bf x}\in F\backslash Q\! :\ g({\bf x})=x_i\}.$$
Of course, then $\pi_{i}(F'_i)$ is the projection of $F'_{i}$ onto the nearest coordinate hyperplane. Therefore, we have
$$\sum\limits_{i=1}^{4}{\rm vol}(D_{i})\geq \iiint\limits_{F\backslash Q} \lambda g({\bf x})d{\bf x}.$$
Let $g_1({\bf x})$ denote the distance from ${\bf x}$ to the nearest face of $F$ (see Fig. \ref{region-center} (1)), then
$$\iiint\limits_{F\backslash Q} \lambda g({\bf x})d{\bf x}=\iiint\limits_{F\backslash Q} \lambda_1 g_1({\bf x})d{\bf x}$$
holds with a suitable constant $\lambda_1$. Obviously, the four sets $F'_1$, $F'_2$, $F'_3$ and $F'_4$ are subsets of the four pyramids with the center ${\bf c}$ of $F$ as their common vertices and with the four facets of $F$ as their bases, respectively.
\begin{figure}
\includegraphics[height=4.5cm]{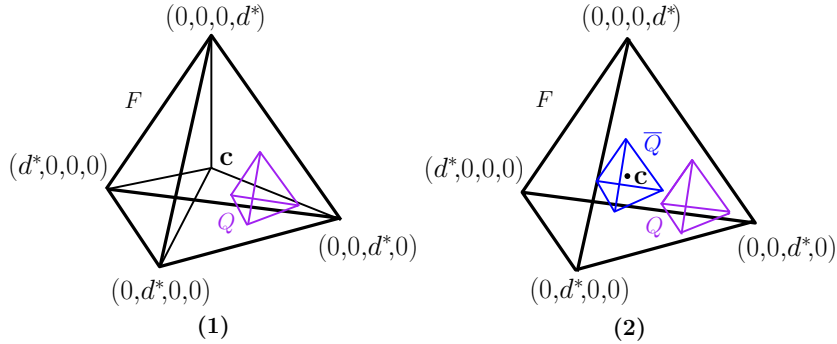}
\caption{The $F'_i$ regions and the placement of $Q$}
\label{region-center}
\end{figure}
Now, we proceed to minimize the integral
$$I=\iiint\limits_{F\backslash Q} \lambda_1g_1({\bf x})d{\bf x}.$$

For convenience, we write
$$h=g_1({\bf c})$$
and define
$$ \chi ({\bf x}, t)=\left\{
\begin{array}{rcl}
1, & & {g_1({\bf x})\geq t},\\
0, & & {g_1({\bf x})< t}.
\end{array}\right. $$
Then, we have
$$I=\lambda_1\iiint\limits_{F\backslash Q}\int_0^h \chi ({\bf x}, t) dtd{\bf x}
=\lambda_1\int_0^h\iiint\limits_{F\backslash Q} \chi ({\bf x}, t)d{\bf x}dt.$$
It is easy to see, for any $t\in [0, h]$, the set $$T=\{{\bf x}\! :\ g_1({\bf x})\geq t\}$$ is a tetrahedron centered at ${\bf c}$, i.e., the blue tetrahedron in Fig. \ref{region-center} (2). Thus, we have
$$\iiint\limits_{F\backslash Q} \chi ({\bf x}, t) d{\bf x}={\rm vol}((F\backslash Q)\cap T).$$
Since the volume of $Q$ is fixed, the integral is minimized when $Q$ is centered at ${\bf c}$. In other words,
$$I=\lambda_1\int_0^h\iiint\limits_{F\backslash Q} \chi ({\bf x}, t) d{\bf x} dt\geq \lambda_1\int_0^h\iiint\limits_{ F\backslash\overline{Q}} \chi ({\bf x}, t) d{\bf x}dt,$$
where $\overline{Q}$ is the translate of $Q$ centered at ${\bf c}$.

Now the only question remaining is, when $Q$ is centered at ${\bf c}$, how large should $Q$ be in order to maximize
$$\frac{d^{*4}}{24}-I-{\rm vol}(P). $$
Suppose that the notch ${\bf v}=(v, v, v, v)$, where $0\leq v\leq \frac{d^*}{4}$, then we have
$${\rm vol}(P)=\frac{(d^*-4v)^4}{4!}.$$
The integral region $\pi_{4}(F'_4)$ and its projection on the $x_1x_2$-plane are shown in Fig. \ref{4-2}. Then, one can deduce that
\begin{align}
\iiint\limits_{\pi_4(F'_4)}\! \left(d^*\!-\!x_1\!-\!x_2\!-\!x_3\right)\!d{\bf x} & =\! \int_0^v \!dx_1\!\int_{d^*-3x_1}^{d^*-x_1}\!dx_2\!\int_{\frac{1}{2}(d^*-x_1-x_2)}^{d^*-x_1-x_2}\! \left(d^*\!-\!x_{1}\!-\!x_{2}\!-\!x_{3}\right)\!dx_3 \nonumber\\
&\quad\! +\!\int_{v}^{d^*\!-3v}\!dx_1\!\int_{d^*\!-x_1-2v}^{d^*\!-x_1}\!dx_{2}\!\int_{\frac{1}{2}(d^*\!-x_1-x_2)}^{d^*\!-x_1-x_2} \!\left(d^*\!-\!x_{1}\!-\!x_{2}\!-\!x_{3}\right)\! dx_{3} \nonumber\\
&\quad\! +\!\int_{d^*\!-3v}^{d^*}\!dx_1\!\int_{\frac{1}{3}(d^*\!-x_1)}^{d^*\!-x_1}\!dx_2\!\int_{\frac{1}{2}(d^*\!-x_1-x_2)}^{d^*\!-x_1-x_2}\!\left( d^*\!-\!x_1\!-\!x_2\!-\!x_3\right)\! dx_3 \nonumber\\
&\quad\! +\!\int_v^{d^*\!-3v}\!dx_1\!\int_v^{d^*\!-x_1-2v}\!dx_2\!\int_{d^*\!-x_1-x_2-v}^{d^*\!-x_1-x_2}\!\left( d^*\!-\!x_1\!-\!x_2\!-\!x_3\right)\! dx_3 \nonumber\\
&\quad\! +\!\int_0^v\!dx_1\!\int_{x_1}^{d^*\!-3x_1}\!dx_2\!\int_{d^*\!-2x_1-x_2}^{d^*\!-x_1-x_2}\!\left( d^*\!-\!x_1\!-\!x_2\!-\!x_3\right)\! dx_3 \nonumber\\
&\quad\! +\!\int_0^v\!dx_2\!\int_{x_2}^{d^*\!-3x_2}\!dx_1\!\int_{d^*\!-x_1-2x_2}^{d^*\!-x_1-x_2}\!\left( d^*\!-\!x_1\!-\!x_2\!-\!x_3\right)\! dx_3 \nonumber\\
&=2v^4-\frac{4d^*}{3}v^3+\frac{d^{*2}}{4}v^2\nonumber.
\end{align}

\begin{figure}
\centering
\includegraphics[height=5cm]{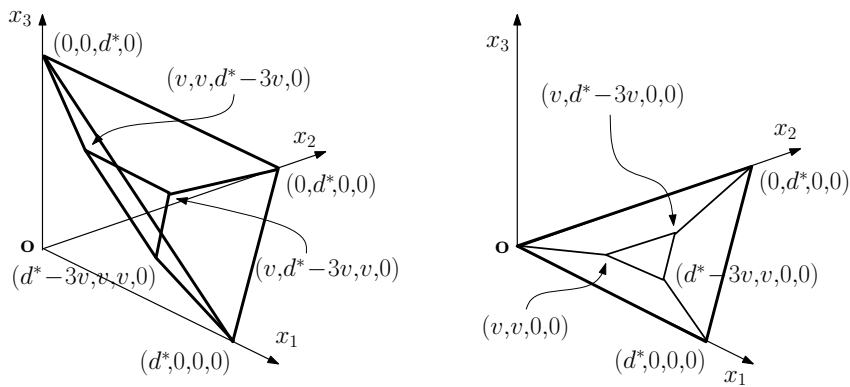}
\caption{The integral region $\pi_{4}(F'_4)$ and its projection on the $x_1x_2$-coordinate plane when there is a notch}
\label{4-2}
\end{figure}

Therefore, in the case that $\overline{T_4}$ has a notch, we have shown that
$${\rm vol}(\overline{T_4})\leq-\frac{56}{3}v^4+16d^*v^3-5d^{*2}v^2+\frac{2d^{*3}}{3}v.\eqno(2.9) $$
To find the maximum of the right-hand side of (2.9), taking its derivative and setting it equal to zero, we get
$$\left(v-\frac{d^*}{4}\right)^2\left(v-\frac{d^*}{7}\right)=0.$$
It is easy to check that the right-hand side of (2.9) attains a local minimum $\frac{d^{*4}}{32}$ at $v=\frac{d^*}{4}$ and attains its maximum  $\frac{11d^{*4}}{343}$ at $v=\frac{d^*}{7}$. In other words, we have
$${\rm vol}(\overline{T_4})\leq \frac{11d^{*4}}{343}\le \frac{11(d+4)^4}{343}.$$

As a conclusion of these two cases, by (2.2) and (2.3) we have shown that
$$f(4,d)\leq \frac{11(d+4)^4}{343}.$$
Theorem 1.3 is proved.

\vspace{1cm}
\noindent
{\LARGE\bf 3. Lattice Coverings of Simplices}

\medskip\noindent
If ${\bf a}_1$, ${\bf a}_2$, $\ldots$, ${\bf a}_n$ are $n$ independent vectors in $\mathbb{E}^n$, then the discrete set
$$\Lambda=\left\{\sum z_i {\bf a}_i\! : z_i\in \mathbb{Z}\right\}$$
is called an $n$-dimensional lattice.  Usually, $\{{\bf a}_1$, ${\bf a}_2$, $\ldots$, ${\bf a}_n\}$ is called a basis of $\Lambda$. Assume that $A$ is the $n\times n$ matrix whose $i$th row is the coordinates of ${\bf a}_i$, then $|{\rm det} A|$ is called the {\it determinant} of $\Lambda$. Usually, it is written as $d(\Lambda)$. If ${\bf a}_1$, ${\bf a}_2$, $\ldots$, ${\bf a}_n$ are $n$ independent vectors in $\mathbb{Z}^n$, as shown before, then $\Lambda$ is an $n$-dimensional sublattice of $\mathbb{Z}^n$. In addition, we have
$$d(\Lambda )=\left[\mathbb{Z}^n\! :\Lambda \right].$$

For a convex body $K$ in $\mathbb{E}^n$, we call $K+\Lambda$ a {\it lattice covering} of $\mathbb{E}^n$ if $$\mathbb{E}^n=\bigcup\limits_{{\bf v}\in \Lambda}(K+{\bf v}).$$
Usually, the value
$$\theta(K, \Lambda)=\frac{{\rm vol}(K)}{d(\Lambda)}$$
is called the density of the covering. Let $\mathcal{L}$ denote the family of all lattices $\Lambda$ such that $K+\Lambda$ is a lattice covering of $\mathbb{E}^n$. Then we call
$$\theta^l(K)=\underset{\Lambda \in \mathcal{L}}{\min}\ \theta (K, \Lambda)=\underset{\Lambda \in \mathcal{L}}{\min} \frac{{\rm vol}(K)}{d(\Lambda)} $$
the {\it lattice covering density} of $K$. Clearly, it is the density of the thinnest lattice covering of $K$.

Given a positive integers $n$ and $d$, we define
$$S^\circ_{n,d}=\{(z_1,\ldots, z_n)\in \mathbb{Z}^n\! :\ z_i\geq 0,\ z_1+\ldots +z_n\leq d\}. $$
Note that the greatest Manhattan distance between a point of $S^\circ_{n,d}$ and ${\bf o}$ is $d$, and
$$\left|S^\circ_{n,d}\right|=\dbinom{d+n}{n}=\frac{d^n}{n!}+O(d^{n-1}).\eqno(3.1)$$
Usually, we call $S^\circ_{n,d}+L_n$ a {\it lattice covering} of $\mathbb{Z}^n$ if $L_n$ is a lattice and
$$\mathbb{Z}^n=\bigcup\limits_{{\bf v}\in L_n}\left(S^\circ_{n,d}+{\bf v}\right).$$
For such a lattice covering we define
$$\theta^*(S^\circ_{n,d}, L_n)=\frac{|S^\circ_{n,d}|}{d(L_n)}.$$
Let $\mathcal{L}_n^*$ denote the family of all $L_n$ such that $S^\circ_{n,d}+L_n$ is a lattice covering of $\mathbb{Z}^n$. Then we call
$$\theta^*(S^\circ_{n,d})=\underset{L_n \in \mathcal{L}_n^*}{\min} \theta^*(S^\circ_{n,d}, L_n)=\underset{L_n \in \mathcal{L}_n^*}{\min} \frac{|S^\circ_{n,d}|}{d(L_n)} $$
the {\it lattice covering density} of $S^\circ_{n,d}$.

Clearly, $S_{n,d}$ is the convex hull of $S^\circ_{n,d}$. To show a useful connection between $\theta^*(S^\circ_{n,d})$ and $\theta^l(S_n)$, we need the following lemma:

\medskip\noindent
{\bf Lemma 3.1 (Dougherty and Faber \cite{Doug04}).} {\it Let $L_n$ be a sublattice of $\mathbb{Z}^n$ and let $\Lambda$ be a lattice in $\mathbb{E}^n$ which is generated by $\{{\bf a}_1$, ${\bf a}_2$, $\ldots$, ${\bf a}_n\}$. Then the following holds:

\begin{enumerate}
\item If $S^\circ_{n,d}+L_n=\mathbb{Z}^n$, then $S_{n, d+n}+L_n=\mathbb{E}^n$.
\item If $S_{n,d}+\Lambda=\mathbb{E}^n$, then there is a constant $c$ such that for all sufficiently large real numbers $k$, if ${\bf b}_i$ is obtained from $k {\bf a}_i$ by rounding all coordinates to the nearest integer, and $L_n$ is the lattice generated by $\{{\bf b}_1$, ${\bf b}_2$, $\ldots$, ${\bf b}_n\}$, then
    $$S_{n,kd+c}+L_n=\mathbb{E}^n.$$
\end{enumerate}}

\medskip\noindent
{\bf Lemma 3.2.} {\it For fixed $n$ and large $d$, we have}
$$\theta^*(S^\circ_{n,d})=\theta^l(S_n)+O(d^{-1}).$$

\smallskip\noindent
{\bf Proof.} Suppose that $S^\circ_{n,d}+L_n$ is a lattice covering of $\mathbb{Z}^n$ with density $\theta^*(S^\circ_{n,d})$. By Lemma 3.1 (1), we have $$S_{n,d+n}+L_n=\mathbb{E}^n.$$
Then, by (3.1) we get
$$\theta^l(S_n)=\theta^l(S_{n,d+n})\leq \theta(S_{n,d+n}, L_n)=\frac{{\rm vol}(S_{n,d+n})\theta^*(S^\circ_{n,d})}{|S^\circ_{n,d}|}=\frac{\theta^*(S^\circ_{n,d})}{1+O(d^{-1})},$$
which implies
$$\theta^*(S^\circ_{n,d})\geq \theta^l(S_n)+O(d^{-1}).\eqno(3.2)$$

On the other hand, suppose that $S_{n,1}+\Lambda$ is a lattice covering of $\mathbb{E}^n$ with density $\theta^l(S_n)$. Let $\{{\bf a}_1$, ${\bf a}_2$, $\ldots$, ${\bf a}_n\}$ and $c$ as chosen in Lemma 3.1. Assume that $d$ is a large integer, we take $k=d-c$ and define $L_n$ to be the lattice approximating $k\Lambda$ as in Lemma 3.1 (2), generated by $\{ {\bf b}_1$, ${\bf b}_2$, $\ldots$, ${\bf b}_n\}$. Now $L_n$ is a sublattice of $\mathbb{Z}^n$ satisfying
$$S_{n,k+c}+L_n=\mathbb{E}^n,$$
so we have
$$S^\circ_{n,d}+L_n=\mathbb{Z}^n.$$
Let $A$ and $B$ be the $n\times n$ matrices with rows ${\bf a}_i$ and ${\bf b}_i$, respectively. Then, we have
$$\theta^*(S^\circ_{n,d})\leq \theta^* (S^\circ_{n,d}, L_n)=\frac{|S^\circ_{n,d}|}{|{\rm det}B|},$$
which together with (3.1) implies
$$\theta^*(S^\circ_{n,d}) \leq \left(\frac{1}{n!}+O(d^{-1})\right)\frac{1}{|{\rm det}(d^{-1}B)|}. \eqno(3.3)$$
By the definitions of $A=(a_{i,j})$ and $B=(b_{i,j})$, we have
$$b_{i,j}=ka_{i,j}+O(1)=da_{i,j}+O(1),\quad i,j =1, 2, \ldots, n,$$
which implies that
$$d^{-1}b_{i,j}=a_{i,j}+O(d^{-1}), \quad i,j =1, 2, \ldots, n.$$
Combined with the assumption that
$$\theta^l(S_n)=\frac{{\rm vol}(S_{n,1})}{|{\rm det}A|}=\frac{1}{n!|{\rm det}A|},$$ we get
$$|{\rm det}(d^{-1}B)|=|{\rm det}A|+O(d^{-1})=\frac{1}{n!\theta^l(S_n)}+O(d^{-1}).\eqno(3.4)$$
Therefore, by (3.3) and (3.4) we get
$$\theta^*(S^\circ_{n,d})\leq \theta^l(S_n)+O(d^{-1}). \eqno(3.5)$$

Clearly, (3.2) and (3.5) together implies
$$\theta^*(S^\circ_{n,d})= \theta^l(S_n)+O(d^{-1}). $$
The lemma is proved. \hfill{$\Box$}

\medskip\noindent
{\bf Remark 3.1.} According to Lemma 3.2, we know that
$$\theta^l(S_n)=\theta^l(S_{n,1})=\lim_{d\to\infty}\theta^*(S^\circ_{n,d}),$$
which is more convenient for us to study $\theta^l(S_n)$.

\medskip
Recall that $T_n$ was defined by (2.1) with M-diameter $d$. To show Theorem 1.1 and Theorem 1.2, we need another result.

\medskip\noindent
{\bf Lemma 3.3.} {\it $S^\circ_{n,d}+L_n$ is a covering of $\mathbb{Z}^n$ if and only if $T_n+L_n$ is a tiling of $\mathbb{Z}^n$, where}
$$T_n=S^\circ_{n,d}\setminus \bigcup_{{\bf v}\in (L_n\cap\mathcal{O}_n)\setminus \{{\bf o}\}}S^\circ_{n,d}+{\bf v}.$$

\medskip\noindent
{\bf Proof.} Since $T_n\subseteq S^\circ_{n,d}$, the sufficiency is obvious. Now we proceed to show the necessary part. In other words, if $S^\circ_{n,d}+L_n$ is a covering of $\mathbb{Z}^n$, for every point $\mathbf{x}\in \mathbb{Z}^n$ there is exact one point ${\bf u}\in L_n$ satisfying ${\bf x}\in T_n+{\bf u}$.

Without loss of generality, we only deal with the points ${\bf x}=(x_1, x_2, \ldots, x_n)\in \mathbb{Z}^n$ with $x_i\ge d$ for all $i=1,$ $2,$ $\ldots,$ $n.$ Let ${\bf u}$ be the $\prec$-last point of $L_n$ such that ${\bf x}\in S^\circ_{n,d}+{\bf u}$ and write
$${\bf y}={\bf x}-{\bf u}\in S^\circ_{n,d}.$$
Then, for any point ${\bf v}\in (L_n\cap\mathcal{O}_n)\setminus \{{\bf o}\}$, we have
$${\bf y}\notin S^\circ_{n,d}+{\bf v}.$$
Otherwise, we obtain that
$${\bf x}={\bf y}+{\bf u}\in S^\circ_{n,d}+{\bf v}+{\bf u},$$
where ${\bf u}\prec {\bf v}+{\bf u}$, which contradicts the $\prec$ assumption on ${\bf u}$. Thus, we must have ${\bf y}\in T_n$ and
$${\bf x}\in T_n+{\bf u}.$$

Suppose that
$${\bf x}={\bf y}+{\bf u}={\bf y}'+{\bf u}',$$
where ${\bf u}$ and ${\bf u}'$ are distinct points of $L_n$ and both ${\bf y}$ and ${\bf y}'$ are in $T_n$ (also in $S^\circ_{n,d}$).
Then, one can deduce that
$${\bf v}={\bf y}-{\bf y}'={\bf u}'-{\bf u}$$
is a nonzero point of $L_n$. Therefore, we have either ${\bf o}\prec {\bf v}$ or ${\bf v}\prec {\bf o}$. In the former case, we have
$${\bf y}={\bf y}'+{\bf v}\in S^\circ_{n,d}+{\bf v},$$
which contradicts to ${\bf y}\in T_n$; In the latter case, we get
$${\bf y}'={\bf y}-{\bf v}\in S^\circ_{n,d}+(-{\bf v}),$$
which contradicts to ${\bf y}'\in T_n$.
The lemma is proved. \hfill{$\Box$}

\medskip\noindent
{\bf Proof of Theorem 1.1.} Based on Proposition 2.1 and Lemma 3.3, we know that the abelian Cayley digraph of $G$ and $\{{\bf g}_1$, ${\bf g}_2$, ${\bf g}_3\}$ has a diameter of at most $d$ if and only if $S^\circ_{3,d}+L_3=\mathbb{Z}^3$. By (1.2) and (2.2), we have
$$\left[\mathbb{Z}^3\! :L_3\right]\leq \frac{3(d+3)^3}{25},$$
which implies
$$\theta^*(S^\circ_{3,d})=\underset{L_3\in \mathcal{L}_3^*}{\min} \frac{|S^\circ_{3,d}|}{\left[\mathbb{Z}^3\! :L_3\right]}\geq\frac{25\binom{d+3}{3}}{3(d+3)^3},$$
where $\mathcal{L}_3^*$ is the family of all $L_3$ such that $S^\circ_{3,d}+L_3$ is a lattice covering of $\mathbb{Z}^3$.
By Remark 3.1, we have
$$\theta^l(S_3)=\lim_{d\to\infty}\theta^*(S^\circ_{3,d})\geq \lim_{d\to\infty}\frac{25\binom{d+3}{3}}{3(d+3)^3}=\frac{25}{18}. $$
Theorem 1.1 is proved. \hfill{$\Box$}

\medskip\noindent
{\bf Proof of Theorem 1.2.} Similar to Theorem 1.1, by (2.2) and Theorem 1.3, we have
$$\left[\mathbb{Z}^4\! :L_4\right]\leq \frac{11(d+4)^4}{343},$$
which implies
$$\theta^*(S^\circ_{4,d})=\underset{L_4\in \mathcal{L}_4^*}{\min} \frac{|S^\circ_{4,d}|}{\left[\mathbb{Z}^4\! :L_4\right]}\geq \frac{343\binom{d+4}{4}}{11(d+4)^4}, $$
where $\mathcal{L}_4^*$ is the family of all $L_4$ such that $S^\circ_{4,d}+L_4$ is a lattice covering of $\mathbb{Z}^4$.
By Remark 3.1, we have
$$\theta^l(S_4)=\lim_{d\to\infty}\theta^*(S^\circ_{4,d})\geq \lim_{d\to\infty}\frac{343\binom{d+4}{4}}{11(d+4)^4}=\frac{343}{264}. $$
Theorem 1.2 is proved. \hfill{$\Box$}

\vspace{0.6cm}\noindent
{\bf Acknowledgements.} This work is supported by the National Natural Science Foundation of China (NSFC11921001, NSFC11801410 and NSFC11971346), the Natural Key Research and Development Program of China (2018YFA0704701), and Jiangsu Natural Science Foundation (BK20210555).}

\bigskip\medskip
\noindent
Miao Fu, Center for Applied Mathematics, Tianjin University, Tianjin, China

\noindent
Email: miaofu@tju.edu.cn

\medskip\noindent
Fei Xue, School of Mathematical Sciences, Nanjing Normal University, Nanjing, China

\noindent
Email: 05429@njnu.edu.cn

\medskip\noindent
Chuanming Zong, Center for Applied Mathematics, Tianjin University, Tianjin 300072, China

\noindent
Email: cmzong@tju.edu.cn

\end{document}